\documentclass[12pt]{article}  
\usepackage{times,amsmath,amsfonts,amsthm,amssymb,eucal}
\usepackage{array}
\usepackage{mathrsfs}
\usepackage{amsmath}
\usepackage{setspace}
\usepackage{fullpage}
\usepackage{mathtools}
\usepackage{fancyhdr}
\usepackage{setspace}
\usepackage{amssymb}
\usepackage{epsfig}
\usepackage{graphicx}
\usepackage{makeidx}
\usepackage{color}
\usepackage{multirow}
\usepackage{epsfig}
\usepackage{graphicx}
\usepackage{makeidx}
\usepackage{gensymb}
\usepackage{enumerate}
\usepackage{multirow}
 \usepackage{setspace}
\usepackage{amssymb}
 \usepackage{setspace}
\newtheorem{theorem}{Theorem}
\newtheorem{proposition}{Proposition}
\newtheorem{lemma}{Lemma}
\newtheorem{corollary}{Corollary}
\newtheorem{definition}{Definition}
\def\Theorem{\begin{theorem}\sl}
\def\EndTheorem{\end{theorem}}
\def\Proposition{\begin{proposition}\sl}
\def\EndProposition{\end{proposition}}
\def\Lemma{\begin{lemma}\sl}
\def\EndLemma{\end{lemma}}
\def\Corollary{\begin{corollary}\sl}
\def\EndCorollary{\end{corollary}}
\def\Definition{\begin{definition}\sl}
\def\EndDefinition{\end{definition}}
\numberwithin{equation}{section}

\begin{document}
\title{ \textbf{On Some Asymptotic Properties and an Almost Sure Approximation of the Normalized Inverse-Gaussian Process}}   
\author{Luai Al Labadi and Mahmoud Zarepour\thanks{{\em Address for correspondence}: M. Zarepour, Department of Mathematics and Statistics,
University of Ottawa, Ottawa, Ontario, K1N 6N5, Canada. E-mail: zarepour@uottawa.ca.} }        
\date{\today}    
\maketitle
\pagestyle {myheadings} \markboth {} {Al Labadi and Zarepour: Asymptotic Properties of the Normalized Inverse-Gaussian Process}
\begin{abstract}
In this paper,  we present some asymptotic properties of the normalized inverse-Gaussian process. In particular, when the concentration parameter is large, we establish an analogue of the empirical  functional central limit theorem, the strong law of large numbers and the Glivenko-Cantelli theorem  for the normalized inverse-Gaussian process and its corresponding quantile process. We also derive a finite sum-representation that converges almost surely to the Ferguson and Klass representation of the normalized inverse-Gaussian process. This almost sure approximation can be used to simulate efficiently the  normalized inverse-Gaussian process.\par

\vspace{9pt}


\vspace{9pt}

\end{abstract}

\section {Introduction}
\label{intro}
The objective of Bayesian nonparametric inference is to  place a prior on the space of probability measures. The Dirichlet process, formally introduced in Ferguson (1973), is considered the first  celebrated  example on this space. Several alternatives of the Dirichlet process have been proposed in the literature. In this paper, we focus on one such prior, namely the normalized inverse-Gaussian process introduced by Lijoi, Mena and Pr\"unster  (2005).

We begin by recalling the definition of the normalized inverse-Gaussian distribution. The random vector  $(Z_1,\ldots,Z_{m})$ is said to have the {\emph{normalized inverse-Gaussian}} distribution with parameters $(\gamma_1,\ldots,\gamma_{m})$, where $\gamma_i > 0$ for all $i$, if it has the joint probability density function
\begin{eqnarray}
\nonumber f(z_1,\ldots,z_{m}) &=& \frac{e^{\sum_{i=1}^\infty\gamma_i} \prod_{i=1}^m\gamma_i}{2^{m/2-1}\pi^{m/2}}\times K_{-m/2}\left(\sqrt{\sum_{i=1}^m\frac{\gamma^2_i}{z_i}}\right) \times \left(\sum_{i=1}^m\frac{\gamma^2_i}{z_i}\right) ^{-m/4}\\ \label{norm}  && \times   \prod_{i=1}^m z_i^{-3/2}\times I_\mathbb{S}(z_1,\ldots,z_{m})
\end{eqnarray}
where $K$ is the modified Bessel function of the third type and $\mathbb{S}=\big\{ (z_1,\ldots,z_{m}): z_i \ge 0,$ $\sum_{i=1}^{m}z_i=1\big\}.$ For more details about the modified bessel functions consult Abramowitz and Stegun (1972, Chapter 9).

Consider a space $\mathfrak{X}$ with a $\sigma-$algebra
$\mathcal{A}$ of subsets of $\mathfrak{X}$. Let $H$ be a fixed probability measure on $(\mathfrak{X},\mathcal{A})$ and $a$ be a positive number.  Following Lijoi, Mena and Pr\"unster (2005), a random probability measure ${P}_{H,a}=\left\{{P}_{H,a}(A)\right\}_{A \in \mathcal{A}}$ is called a
 normalized inverse-Gaussian process on $(\mathfrak{X},\mathcal{A})$ with parameters
$a$ and $H$, if for any finite measurable partition $A_1, \ldots,$ $
A_m$ of $\mathfrak{X}$, the joint distribution of the vector
$\left({P}_{H,a}(A_1), \ldots\,{P}_{H,a}(A_{m})\right)$ has the normalized inverse-Gaussian distribution with parameter $ \left(a H(A_1), \ldots\,a H(A_m)\right)$. We assume that if $H(A_i)=0$, then ${P}_{H,a}(A_i)=0$ with probability one. The normalized inverse-Gaussian process with parameters $a$ and $H$ is denoted by $\text{N-IGP}(a, H)$, and we write ${P}_{H,a}\sim \text{N-IGP}(a, H).$

One of the basic properties of the normalized inverse-Gaussian process is that for any   $A \in \mathcal{A},$
  \begin{equation}
{E}({P}_{H,a}(A))=H(A) \ \ \text{ and }\ \ {Var}({P}_{H,a}(A))=\frac{H(A)(1-H(A))}{\xi(a)}, \label{NIG1}
\end{equation}
where here and throughout this paper
  \begin{equation}
  \xi(a)=\frac{1}{a^2 e^a \Gamma(-2,a)} \label{xi}
  \end{equation}
and $\Gamma(-2,\theta)=\int_{a}^\infty t^{-3}e^{-t}dt$.

Furthermore,  for any two disjoint sets $A_i$ and $A_j \in \mathcal{A},$
\begin{equation}
{E}({P}_{H,a}(A_i){P}_{H,a}(A_j))=H(A_i)H(A_j)\frac{\xi(a)-1}{\xi(a)}. \label{As.20}
\end{equation}

Observe that, for large $a$, $\xi(a) \approx {a}$ (Abramowitz and Stegun, 1972, Formula 6.5.32, page 263), where we use the notation $f(a) \approx g(a)$ if $\lim_{a \to \infty} {f(a)}/{g(a)}=1.$ It follows from (\ref{NIG1}) that $H$ plays the role of the center
 of the process, while  $a$  can be viewed as  the concentration parameter. The larger $a$ is, the more likely it is that the realization of ${P}_{H,a}$ is close to $H$. Specifically, for any fixed set  $A \in \mathcal{A}$ and $\epsilon>0,$ we have ${P}_{H,a}(A) \overset{p} \to H(A)$ as $a \to \infty$ since
\begin{equation}
\Pr\left\{\left|{P}_{H,a}(A)-H(A)\right|>\epsilon\right\} \le \frac{H(A)(1-H(A))}{\xi(a)\epsilon^2}.\label{cheb0}
\end{equation}

Similar to the Dirichlet process, a series representation of the normalize inverse-Gaussian process can be easily derived  from the Ferguson and Klass representation (1972). Specifically, let $(E_i)_{i \ge 1}$ be a sequence of independent and identically distributed (i.i.d.) random variables with an exponential distribution with mean of 1. Define
\begin{equation}
\Gamma_i=E_1+\cdots+E_i. \label{eq2}
\end{equation}
Let $(\theta_i)_{i \geq 1}$  be a sequence of i.i.d. random variables with values in $\mathfrak{X}$ and common distribution $H$, independent of $(\Gamma_i)_{i
\ge 1}$. Then  the normalize inverse-Gaussian process with parameter $a$ and $H$ can be expressed  as a normalized series representation
\begin{equation}
P_{H,a}(\cdot)=\sum_{i=1}^{\infty} {\frac{L^{-1}(\Gamma_i)}{\sum_{i=1}^{\infty}{{L^{-1}(\Gamma_i)}}}\delta_{\theta_i}(\cdot)},\label{eq3}
\end{equation}
where
\begin{equation}
L(x)=\frac{a}{\sqrt{2\pi}} \int_{x}^{\infty}{e^{-t/2}}{t^{-3/2}}~dt, \text{ for } x>0, \label{eq5}
\end{equation}
and $\delta_X$ denotes the Dirac measure at $X$ (i.e. $\delta_X(B)=1$
if $X \in B$ and $0$ otherwise). Note that, working with (\ref{eq3}) is difficult in practice because no closed form for the inverse of the L\'evy measure (\ref{eq5}) exists. Moreover, to determine the random weights in (\ref{eq3}) an infinite sum must be computed.

This paper is organized as follows. In Section 2, we study  the weak convergence of the centered and scaled process
\begin{equation}
{D}_{H,a}(\cdot)=\sqrt{\xi(a)}\left({P}_{H,a}(\cdot)-H(\cdot)\right), \nonumber
\end{equation}
as $a \to \infty.$ Note that, since $\xi(a)\approx a$, it is possible to replace the normalizing coefficient $\xi(a)$ by $a$. Therefore, for simplicity, we focus on the process
\begin{equation}
{D}_{H,a}(\cdot)=\sqrt{a}\left({P}_{H,a}(\cdot)-H(\cdot)\right). \label{NIGP2.18}
\end{equation}
In Section 3, we derive the limiting process  for the  quantile process
\begin{equation}
{Q}_{H,a}(\cdot)=\sqrt{a}\left({{P}}^{-1}_{H,a}(\cdot)-H^{-1}(\cdot)\right),\label{NIG4}
\end{equation}
as $a \to \infty,$ where, in general, the inverse of a distribution function $F$ is defined by
\begin{equation}
F^{-1}(t)=\inf \left\{x:F(x)\ge t \right\}, \ \ 0<t<1. \nonumber
\end{equation}
The strong law of large numbers and the Glivenko-Cantelli theorem  for the normalized inverse-Gaussian process are discussed in Section 4. In Section 5, we derive a finite sum-representation which converges almost surely to the Ferguson and Klass representation of the normalized inverse-Gaussian process. This new representation provides a simple, yet efficient procedure, to sample the the normalized inverse-Gaussian process.

\section{Asymptotic Properties of the Normalized Inverse-Gaussian Process}
In this section,  we study the weak convergence of  the process ${D}_{H,a}$  defined in  (\ref{NIGP2.18}) for large values of $a$. Let $\mathscr{S}$  be a collection of Borel sets in $\mathbb{R}$  and $H$ be a probability measure on $\mathbb{R}$. We  recall the definition of a  Brownian bridge indexed by $\mathscr{S}.$ A Gaussian process $\left\{B_H(S): S \in \mathscr{S}\right\}$ is called a {\emph{Brownian bridge with parameter measure }$H$  if ${E}\left[B_H(S)\right]=0$ for any $S \in \mathscr{S}$ and
\begin{eqnarray}
\label{asym:cov0}  {Cov}\left(B_H(S_i),B_H(S_j)\right)=H(S_i \cap S_j)- H(S_i)H(S_j)
 \end{eqnarray}
for any  $S_i, S_j \in \mathscr{S}$ (Kim and Bickel, 2003).

The next lemma gives the limiting distribution of the process  (\ref{NIGP2.18}) for any finite Borel set $S_1,\ldots ,S_m \in \mathscr{S}$, as $a \to \infty$. The proof of the lemma for $m=2$ is given in the appendix and can be generalized easily to the case of arbitrary $m.$

\Lemma \label{Asy.L2} Let $D_{H,a}$ be  defined by  (\ref{NIGP2.18}). For any fixed  sets $S_1,\ldots ,S_m$ in $\mathscr{S}$ we have
\begin{equation}
\left(D_{H,a}(S_1),D_{H,a}(S_2),\ldots,D_{H,a}(S_m)\right)\overset {d }\to \left(B_H(S_1),B_H(S_2),\ldots,B_H(S_m)\right), \nonumber
\end{equation}
as $a \to \infty$, where  $B_H$ is the Brownian bridge with parameter $H$.

\EndLemma

The following theorems show that the process  ${D}_{H,a}$  defined by  (\ref{NIGP2.18}) converges to  the process $B_H$ on  $D[-\infty,\infty]$ with respect to the Skorokhod  topology, where $D[-\infty,\infty]$ is the space of cadlag functions (right continuous with left limits) on $[-\infty,\infty]$. Right continuity at $-\infty$ can be achieved by setting $D_{H,a}(-\infty)={D}_{H,a}(-\infty)=0$; the left limit at $\infty$ also equals zero, the natural value of $D_{H,a}(\infty)$ and ${D}_{H,a}(\infty)$.  For more details, consult Pollard (1984, Chapter 5). If $X$ and $(X_a)_{a>0}$ are random variables with values in a metric space
$M$, we say that $(X_a)_{a}$ converges in distribution to $X$ as $a\to \infty$ (and we write $X_{a}
\overset{d}\to X)$ if for any sequence $(a_n)_n$ converging to $\infty$,
$X_{a_n}$ converges in distribution to $X.$

\Theorem
\label{Asy.TT1} We have, as $a\to \infty$,
 \begin{eqnarray}
D_{H,a}(\cdot)=\sqrt{a}\left(P_{H,a}(\cdot)-H(\cdot)\right) \overset{d} \to B_H(\cdot) \label{eq2.188}
 \end{eqnarray}
 on $D[-\infty,\infty]$ with respect to Skorokhod  topology, where $B_H$ is the Brownian bridge with parameter measure $H$.
\EndTheorem

\noindent {\textbf{Remark 1.}} In the proof of Lemma \ref{Asy.L2}, the finite-dimensional convergence is in fact convergence in total variation, which is stronger than convergence in distribution (Billingsley 1999, page 29).

\noindent {\textbf{Remark 2.}}  Lo (1987)  obtained a result similar to that given in  Theorem \ref{Asy.TT1} for the Dirichlet process to establish asymptotic validity of the Bayesian bootstrap. An interesting generalization of Lo (1987)  to  the two-parameter Poisson-Dirichlet process was obtained by James (2008). In both papers, proofs are based on constructing the distributional identities (Proposition 4.1, James, 2008)). Establishing analogous distributional identity for the normalized inverse-Gaussian process does not seem to be trivial.

\noindent {\textbf{Remark 3.}} Sethuraman and Tiwari (1982) studied the convergence and tightness of the Dirichlet process as the parameters are allowed to converge in a certain sense. Analogous results for  the two-parameter Poisson-Dirichlet process (Pitman and Yor, 1997; Ishwaran and James, 2001) and the normalized inverse-Gaussian process $\left({P}_{H,a}\right)_{a>0}$ follows straightforwardly by applying the technique used in the proof of the tightness part (i.e. condition (\ref{tight})) of Theorem \ref{Asy.TT1}.

\section{Asymptotic Properties of the  Normalized Inverse-Gaussian Quantile Process}
Similar to the frequentist asymptotic theory, in this section we establish large sample theory for the normalized inverse-Gaussian quantile process.
\Corollary \label{Asy.} Let $0<p<q<1, $ and  $H$ be a continuous function with positive derivative $h$ on the interval $\left[H^{-1}(p)-\epsilon,H^{-1}(q)+\epsilon\right]$ for some $\epsilon>0.$ Let $\lambda$ be the Lebesgue  measure on $[0,1]$.
Let ${Q}_{H,a}$ be the normalized inverse-Gaussian process defined in  (\ref{NIG4}). As $a\to \infty$, we have
 \begin{eqnarray}
 \nonumber Q_{H,a}(\cdot) \overset{d} \to -\frac{B_\lambda(\cdot)}{h(H^{-1}(\cdot))}=Q(\cdot),
 \end{eqnarray}
 in $D[p,q].$ That is, the limiting process is a Gaussian process with zero-mean and covariance function
 $$Cov\left(Q(s),Q(t)\right)=\frac{\lambda(s\wedge t)-\lambda(s)\lambda(t)}{h(H^{-1}(s))h(H^{-1}(t))},\quad s, t \in \mathbb{R}.$$
\EndCorollary

\proof We only proof part (i). Part (ii) follows similarly. By Theorem \ref{Asy.TT1} (i) the process $\sqrt{a}\left(P_{H,a}-H\right)$ converges in distribution to the process $B_H=B_\lambda(H)=B_\lambda \circ H.$ Notice that almost all sample paths of the limiting process are continuous  on the interval $\left[H^{-1}(p)-\epsilon,H^{-1}(q)+\epsilon\right].$  By Lemma 3.9.23 page 386  of van der Vaart and Wellner (1996),  the inverse map $H \mapsto H^{-1}$ is Hadamard tangentially differentiable at $H$ to the subspace of functions that are continuous on this interval. By the functional delta method (Theorem 3.9.4 page 374 of van der Vaart and Wellner (1996)) we have
\begin{eqnarray}
 \nonumber Q_{H,a}(\cdot) \overset{d} \to -\frac{B_\lambda\circ H\circ H^{-1}(\cdot)}{h(H^{-1}(\cdot))}=-\frac{B_\lambda(\cdot)}{h(H^{-1}(\cdot))}
 \end{eqnarray}
 in $D[p,q]$. This completes the proof of the corollary.
\endproof

 \noindent {\textbf{Remark 4.}} A similar result to Corollary \ref{Asy.} for the two-parameter Poisson-Dirichlet process  can be obtain by applying Theorem 4.1 and Theorem 4.2 of James (2008).

 \noindent {\textbf{Remark 5.}}  Similar to Remark 1 of Bickel and Freedman (1981), if  $H^{-1}(0+)>-\infty$ and $H^{-1}(1)<\infty$ and $h$ is continuous on $[H^{-1}(0+),$ $H^{-1}(1)],$ the conclusion of the corollary holds in  $D\left[H^{-1}(0+),H^{-1}(1)\right].$  For example, if $H$ is a uniform distribution on $[0,1],$ then the convergence holds in $D[0,1].$

\vspace{2mm}
The next example is a direct application of Corollary \ref{Asy.}.

\noindent {\textbf{Example 1.}}   In this example we derive the asymptotic distribution for the median and the interquantile range for the the normalized inverse-Gaussian process. Let $Q^{1}_{H,a}$, $Q^{2}_{H,a}$ and $Q^{3}_{H,a}$  be the first, the second (median) and the third quartiles of $P_{H,a}$ (i.e. $P^{-1}_{H,a}(0.25)=Q^1_{{H,a}}$,  $P^{-1}_{H,a}(0.5)=Q^2_{{H,a}}$ and $P^{-1}_{H,a}(0.75)=Q^{3}_{H,a}$). Let $q_{1}$, $q_2$ and $q_{3}$ be  the first, the second (median) and the third quartiles of $H$. From Corollary \ref{Asy.},  after some simple calculations, the asymptotic distribution of the median and the  interquantile range are given, respectively, by:
\begin{eqnarray}
 \nonumber \sqrt{a}\left(Q^{2}_{H,a}-q_2\right) \overset{d} \to N\left(0,\frac{1}{4h^2(q_2)}\right)
 \end{eqnarray}
and
 \begin{eqnarray}
 \nonumber \sqrt{a}\left(IQR-(q_3-q_1)\right) \overset{d} \to N\left(0,\frac{3}{h^2(q_3)}+\frac{3}{h^2(q_1)}-\frac{2}{h(q_1)h(q_3)}\right),
 \end{eqnarray}
where $h=H^{\prime}$ and  $IQR=Q^{3}_{H,a}-Q^{1}_{H,a}$. Note that, the asymptotic distributions  of the median and the interquantile range for the normalized inverse-Gaussian process coincide with that of  the sample median and the sample interquartile range (DasGupta, 2008, page 93).

\section {Glivenko-Cantelli Theorem  for the  Normalized Inverse-Gaussian Process}

In this section, we  show that an analogue of the empirical strong law of large numbers and the empirical Glivenko-Cantelli theorem continue to hold for  the normalized inverse-Gaussian process.

\Theorem \label{gliv1} Let ${P}_{H,a}\sim \text{N-IGP}(a,H)$.  Assume that $a=n^2c$, for a fixed positive number $c$.  Then as $n \to \infty,$
 \begin{eqnarray*}
 \label{gliv} {P}_{H,n^2c}(A) \overset{a.s.} \to H(A),
 \end{eqnarray*}
 for any  measurable subset $A$  of $\mathfrak{X}$.

\EndTheorem
\proof  For any $\epsilon>0$, by (\ref{cheb0}), we have
\begin{equation}
\Pr\left\{\left|{P}_{H, n^2c}(A)-H(A)\right|>\epsilon\right\} \le \frac{H(A)(1-H(A))\xi(n^2c)}{\epsilon^2},\label{cheb0}
\end{equation}
where $\xi(n^2c)$ is defined by (\ref{xi}). Note that $$\lim_{n\to \infty}\frac{\xi(n^2c)}{1/n^2c}=1$$
(Abramowitz and Stegun, 1972, Formula 6.5.32, page 263). Since the series $\sum_{n=1}^\infty 1/n^2$ converges, it follows by the Limit Comparison Test that the series  $\sum_{n=1}^\infty \xi(n^2c)$ is also convergent. Thus,
\begin{eqnarray*}
\sum_{n=1}^{\infty} \Pr \left\{|{P}_{H, n^2c}(A)-H(A)|>\epsilon\right\}<\infty.
 \end{eqnarray*}
Therefore, by the first Borel-Cantelli Lemma, the proof follows.
\endproof

The proof of the next theorem follows by arguments similar to that given in the proof of  the Glivenko-Cantelli theorem for the empirical process. See, for example, Billingsely (1995, Theorem 20.6).

\Theorem \label{Asy.T3} Let ${P}_{H,a}\sim \text{N-IGP}(a,H)$.  Assume that $a=n^2c$, for a fixed positive number $c$. Then
 \begin{eqnarray}
 \label{gliv} \sup_{x \in \mathbb{R}}\left|{P}_{H,n^2c}(x)-H(x)\right| \overset{a.s.} \to 0,
 \end{eqnarray}
 as $n \to  \infty$.

 \noindent {\textbf{Remark 6.}}  Similar to the normalized inverse-Gaussian process, a strong law of large numbers and a Glivenko-Cantelli theorem can also be established for the two-parameter Poisson-Dirichlet process.

\EndTheorem

\section{Monotonically Decreasing Approximation to the Inverse Gaussian Process} \label{new}
In this section, we derive a finite sum representation which converges almost surely to the Ferguson and Klass sum representation of the normalized inverse-Gaussian process. We mimic the approach developed recently by Zarepour and Al Labadi (2012) for the Dirichlet process. Let $X_n$  be a random variable with distribution inverse-Gaussian with parameter $a/n$ and $1$ (see equation (3) of Lijoi, Mena and Pr\"unster (2005) for the density of the inverse-Gaussian distribution). Define
\begin{equation}
G_n(x)=\Pr(X_n>x)=\int_{x}^{\infty}\frac{a}{n\sqrt{2\pi}}t^{-3/2}\exp\left\{-\frac{1}{2}\left(\frac{a^2}{n^2t}+t\right)+\frac{a}{n}\right\}dt. \label{eq9}
\end{equation}

The following proposition describes properties of $G_n(x)$ that  will be used later in the paper.

\Proposition \label{prop1}
For $x>0$, the function $G_n(x)$ defined in (\ref{eq9}) has the following properties as $n \to \infty$:
\begin{enumerate}[(i)]
\item $nG_n(x)\rightarrow L(x),$
\item $G_n^{-1}(\frac{x}{n})  \rightarrow L^{-1}(x),$
\end{enumerate}
where $L$ is defined in (2.3).
\EndProposition

\proof
To prove (i), since $e^{a/n}\le e^a$ and $\exp\left\{-\frac{1}{2}\left(\frac{a^2}{n^2t}+t\right)\right\}\le 1$, the integrand in $nG_n(x)$ is bounded by $ae^at^{-3/2}/\sqrt{2\pi}$, which is integrable for any $x>0$. Hence, the dominated convergence theorem applies and we have
\begin{equation}
nG_n(x)=\int_{x}^{\infty}\frac{a}{\sqrt{2\pi}}t^{-3/2}e^{-\frac{1}{2}\left(\frac{a^2}{n^2t}+t\right)+\frac{a}{n}}dt \rightarrow \frac{a}{\sqrt{2\pi}} \int_{x}^{\infty}{e^{-t/2}}{t^{-3/2}}~dt=L(x).\nonumber
\end{equation}

To prove (ii), notice that the left hand side of (i) is a sequence of
monotone functions converging to a continuous monotone function  for every $x>0$ (Haan-de and Ferreira 2006, page 5). Thus, (i) is
equivalent to $G^{-1}_n(x/n)\rightarrow N^{-1}(x).$
\endproof

Proposition \ref{prop1} gives a simple procedure for an approximate evaluation
of both $L(x)$ and $L^{-1}(x)$ for any $x>0$. For computational
simplicity, a more convenient approximation is presented in the
following Corollary.  The proof follows straightforwardly by taking
$x=\Gamma_i$ in Proposition 1 and the fact that we have ${\Gamma_{n+1}}/{n} \overset{a.s.}\rightarrow 1$  as $n \to \infty$ (strong law of large numbers).

\Corollary \label{cor2}
For a fixed $i$, as $n \to \infty,$ we have:
$$L_n^{-1}\left(\frac{\Gamma_i}{\Gamma_{n+1}}\right) \overset{a.s.} \rightarrow L^{-1}(\Gamma_i).$$
\EndCorollary

 \noindent {\textbf{Remark 7.}}  The utility of Corollary \ref{cor2} stems from the fact that all values of
$G_n^{-1}({\Gamma_i}/{\Gamma_{n+1}})$ are nonzero for $i \le
n$. This is not the case when working with
$G_n^{-1}({\Gamma_i}/{n})$.

\vspace{5pt}

The following lemma provides a finite sum representation which converges,
almost surely, to the Ferguson and Klass~(1972) sum-representation for the
inverse-Gaussian process. The proof of the lemma  is similar to that of Lemma 2 in Zarepour and Al Labadi (2012). Hence, it is omitted. 

\Lemma \label{BG1}
If $(\theta_i)_{i \geq 1}$ is a sequence of i.i.d. random variables
with common distribution $H$, independent of $(\Gamma_i)_{i \ge 1}$,
then as $n \to \infty$
\begin{equation}
\sum_{i=1}^{n} {{G_n^{-1}\left(\frac{\Gamma_i}{\Gamma_{n+1}}\right)}\delta_{\theta_i}}\overset{a.s.}\rightarrow
\sum_{i=1}^{\infty} {L^{-1}\left(\Gamma_i\right)\delta_{\theta_i}}. \label{INGN}
\end{equation}
 Here, $\Gamma_i$, $N(x)$, and $G_n(x)$, are defined in (\ref{eq2}),
(\ref{eq5}), and (\ref{eq9}), respectively.
\EndLemma


By normalizing the finite sum in (\ref{INGN}), it is possible to obtain a sum representation that converges almost surely
to Ferguson and Klass representation of the normalized inverse-Gaussian process. This important result
 is stated formally in the next theorem.

\Theorem \label{Theorem2}
Let $(\theta_i)_{i \geq 1}$ be a sequence of i.i.d. random variables
with values in $\mathfrak{X}$ and common distribution $H$, independent
of $(\Gamma_i)_{i \ge 1}$, then as $n \to \infty$
\begin{equation}
P^{\text{new}}_{n,H,a}=\sum_{i=1}^{n} \frac
{{G_n^{-1}\left(\frac{\Gamma_i}{\Gamma_{n+1}}\right)}}{\sum_{i=1}^{n}{G_n^{-1}\left(\frac{\Gamma_i}{\Gamma_{n+1}}\right)}}\delta_{\theta_i}\overset{a.s.}\rightarrow
P_{H,a}=\sum_{i=1}^{\infty} {\frac{L^{-1}(\Gamma_i)}{\sum_{i=1}^{\infty}{{L^{-1}(\Gamma_i)}}}\delta_{\theta_i}}.
\label{eq11}
\end{equation}
Here $\Gamma_i$, $L(x)$,  and $G_n(x)$, are defined in (\ref{eq2}),
(\ref{eq5}), and (\ref{eq9}), respectively.
\EndTheorem

\noindent {\textbf{Remark 7.}}  For any $1 \le i \le n,$  ${\Gamma_{i}}/{\Gamma_{n+1}}<{\Gamma_{i+1}}/{\Gamma_{n+1}}$  almost surely. Since $G_n^{-1}$ is a decreasing function, we have $G_n^{-1}\left({\Gamma_{i}}/{\Gamma_{n+1}}\right)>G_n^{-1}\left( {\Gamma_{i+1}}/{\Gamma_{n+1}}\right)$ almost surely. That is, the weights of the new representation given in Theorem \ref{Theorem2}  decrease monotonically for any fixed positive integer $n$. As demonstrated in Zarepour and Al Labadi (2012) in the case of the Dirichlet process, we anticipate that this new representation will yield highly  accurate approximations to the normalized inverse-Gaussian process.

\noindent {\textbf{Remark 7.}}  For $P^{\text{new}}_{n,H,a}$ of Theorem \ref{Theorem2} we can write
\begin{equation}
P^{\text{new}}_{n,H,a}\overset{d}=\sum_{i=1}^{n} p_{i,n}\delta_{\theta_i},
\label{eq11}
\end{equation}
where $p_{1,n},\ldots,p_{n,n}\sim\text{N-IG}\left(a/n,\ldots,a/n\right)$, $\overset{d}=$ means have the same distribution and the N-IG distribution is given by (\ref{norm}). Therefore a similar result to Theorem 2 of Ishwaran and Zarepour for the normalized inverse-Gaussian process follows immediately.
\vspace{5pt}

\section {Concluding Remarks}
The approach used in this paper  can be applied to similar processes with tractable finite dimensional distributions. On the other hand, when the finite dimensional distribution is unknown, one may follow the approach of James (2008). However, applying this approach requires constructing distributional identities, which may be difficult in some cases.

One can use the results obtained  in this paper to  derive asymptotic properties of any  Hadamard-differentiable functional of the $\text{N-IGP}(a,H)$ as $a\to \infty.$ For different applications in statistics we refer the reader to van der Vaart and Wellner (1996, Section 3.9) and  Lo (1987). Moreover, it is possible to extend the results found in this paper to the case when the base measure $H$ is a multivariate cumulative distribution function.  The result of Bickle and Wichura (1972) can be employed in the proof.

\section{Acknowledgments} The authors would like to thank Professor Raluca Balan for her helpful comments and suggestions. The research of the authors are supported by research grants from the \textbf{Natural Sciences and
Engineering Research Council of Canada (NSERC)}.

\newpage
 \section*{Appendix}
\noindent {\textbf{Proof of Lemma 1 for} $\mathbf{m=2}$\textbf{:}}
Let $S_1$ and $S_2$ be any two intervals in $\mathbb{R}$. Without loss of generality, we assume that $S_1 \cap S_2=\emptyset$. The general case when $S_1$ and $S_2$ are not necessarily disjoint follows from the continuous mapping theorem.

Note that
 \begin{eqnarray}
\nonumber \left({P}_{H,a}(S_1),{P}_{H,a}(S_2),1-{P}_{H,a}(S_1)-{P}_{H,a}(S_2)\right) &\sim &\text{N-IG}\big(a H(S_1),a H(S_2),\\ \nonumber  &&a(1-H(S_1)-H(S_2))\big),
 \end{eqnarray}
 where the N-IG distribution is given by (\ref{norm}).  For notational simplicity, set  $X_i={P}_{H,a}(S_i)$, $l_i=H(S_i)$ and  ${D}_i=\sqrt{a} \left(X_i-l_i\right)$ for $i=1,2.$ Thus, the joint density function of $X_1$ and $X_2$ is:
\begin{eqnarray}
\nonumber f_{X_1,X_2}(x_1,x_2) &=& \frac{e^{a} a^3 l_1 l_2(1-l_1-l_2)}{2^{1/2}\pi^{3/2}}\times x_1^{-3/2} x_2^{-3/2} (1-x_1-x_2)^{-3/2}
\\ \nonumber && \times K_{-3/2}\left(a\sqrt{\frac{l_1^2}{x_1}+\frac{l_2^2}{x_2}+\frac{(1-l_1-l_2)^2}{1-x_1-x_2}}\right)
\\ \nonumber &&  \times a^{-3/2}\left(\frac{l_1^2}{x_1}+\frac{l_2^2}{x_2}+\frac{(1-l_1-l_2)^2}{1-x_1-x_2}\right)^{-3/4}
\\ \nonumber &=&  \frac{a^{1/2}e^{a} l_1 l_2(1-l_1-l_2)}{2^{1/2}\pi^{3/2}}\times  x_1^{-3/2} x_2^{-3/2} (1-x_1-x_2)^{-3/2}
\\ \nonumber && \times K_{-3/2}\left(a\sqrt{\frac{l_1^2}{x_1}+\frac{l_2^2}{x_2}+\frac{(1-l_1-l_2)^2}{1-x_1-x_2}}\right)
\\ \nonumber &&  \times \left(\frac{l_1^2}{x_1}+\frac{l_2^2}{x_2}+\frac{(1-l_1-l_2)^2}{1-x_1-x_2}\right)^{-3/4}.
\end{eqnarray}

The joint probability density function of ${D}_1=\sqrt{a} \left(X_1-l_1\right)$ and  ${D}_2=\sqrt{a} \left(X_2-l_2\right)$ is:
\begin{eqnarray}
\nonumber f_{{D}_1,{D}_2}(y_1,y_2)&=& \frac{a^{1/2}e^{a} l_1 l_2(1-l_1-l_2)}{2^{1/2}\pi^{3/2}}\times
\\ \nonumber && \times \left({y_1}/{\sqrt {a}}+l_1\right)^{-3/2}  \left({y_2}/{\sqrt {a}}+l_2\right)^{-3/2}
\\ \nonumber && \left(1-{y_1}/{\sqrt {a}}-l_1-{y_2}/{\sqrt {a}}-l_2 \right)^{-3/2}
\\ \nonumber && \times K_{-3/2}\Bigg(a\Bigg(\frac{l_1^2}{{y_1}/{\sqrt {a}}+l_1}+\frac{l_2^2}{{y_2}/{\sqrt {a}}+l_2}\\ \nonumber &&+\frac{(1-l_1-l_2)^2}{1-{y_1}/{\sqrt {a}}-l_1-{y_2}/{\sqrt {a}}-l_2}\Bigg)^{1/2}\Bigg)\\
\nonumber &&  \times \Bigg(\frac{l_1^2}{{y_1}/{\sqrt {a}}+l_1}+\frac{l_2^2}{{y_2}/{\sqrt {a}}+l_2}\\
 \nonumber && +\frac{(1-l_1-l_2)^2}{1-{y_1}/{\sqrt {a}}-l_1-{y_2}/{\sqrt {a}}-l_2}\Bigg)^{-3/4}.
\end{eqnarray}
By Scheff\'e's theorem (Billingsely 1999, page 29), it is enough to show that:
\begin{eqnarray}
\label{asy4} f_{D_1,D_2}(y_1,y_2) \to f(y_1,y_2) =\frac{1}{2\pi |\Sigma|^{1/2}}\exp\left\{-(y_1~y_2)\Sigma^{-1}(y_1~y_2)^T/2\right\},
\end{eqnarray}
where $\Sigma=\begin{bmatrix}  l_1\left(1-l_1\right)& -l_1l_2 \\ -l_1l_2 & l_2\left(1-l_2\right) \end{bmatrix}.$\\

Since, for large $z$ and fixed $\nu$, $K_\nu(z) \approx \sqrt{{\pi}/{2}}z^{-1/2}e^{-z}$ (Abramowitz and Stegun, 1972, Formula 9.7.2, page 378), where we use the notation $f(z) \approx g(z)$   if $\lim_{z \to \infty} {f(z)}/{g(z)}=1,$
we get:
\begin{eqnarray}
\nonumber \lim_{a \to \infty}f_{D_1,D_2}(y_1,y_2)&=& \lim_{a \to \infty}\Bigg[\frac{l_1 l_2(1-l_1-l_2)}{2\pi}\times
\\ \nonumber && \Bigg(\frac{l_1^2}{y_1/\sqrt{a}+l_1}+\frac{l_2^2}{y_2/\sqrt{a}+l_2}+\frac{(1-l_1-l_2)^2}{1-y_1/\sqrt{a}-l_1-y_2/\sqrt{a}-l_2}\Bigg)^{-1}
\\ \nonumber && \times \left(y_1/\sqrt{a}+l_1\right)^{-3/2}  \left(y_2/\sqrt{a}+l_2\right)^{-3/2}
\\ \nonumber && \left(1-y_1/\sqrt{a}-l_1- y_2/\sqrt{a}-l_2\right)^{-3/2}
\\ \nonumber && \times \exp\Bigg(a\Bigg(1-\Bigg(\frac{l_1^2}{y_1/\sqrt{a}+l_1}+\frac{l_2^2}{y_2/\sqrt{a}+l_2}
\\\nonumber &&+\frac{(1-l_1-l_2)^2}{1-y_1/\sqrt{a}-l_1-y_2/\sqrt{a}-l_2}\Bigg)^{1/2}\Bigg)\Bigg)\Bigg].
\end{eqnarray}

Notice that,
\begin{eqnarray}
\nonumber &&\frac{l_1 l_2(1-l_1-l_2)}{2\pi}\times \Bigg(\frac{l_1^2}{y_1/\sqrt{a}+l_1}+\frac{l_2^2}{y_2/\sqrt{a}+l_2}+\frac{(1-l_1-l_2)^2}{1-y_1/\sqrt{a}-l_1-y_2/\sqrt{a}-l_2}\Bigg)^{-1}
\\ \nonumber && \times \left(y_1/\sqrt{a}+l_1\right)^{-3/2}  \left(y_2/\sqrt{a}+l_2\right)^{-3/2}\left(1-y_1/\sqrt{a}-l_1- y_2/\sqrt{a}-l_2\right)^{-3/2}
\end{eqnarray}

converges to ${1}/{\left(2\pi \sqrt{\sigma_{11}\sigma_{22}(1-\rho_{12}^2)}\right)},$ where
\begin{eqnarray}
\nonumber \sigma_{11}=l_1(1-l_1),\ \ \sigma_{22}=l_2(1-l_2), \ \  \rho_{12}=-\sqrt{\frac{l_1l_2}{(1-l_1)(1-l_2)}}.
\end{eqnarray}
To prove the lemma, it remains to show that
$$a\Bigg(1-\Bigg(\frac{l_1^2}{y_1/\sqrt{a}+l_1}+\frac{l_2^2}{y_2/\sqrt{a}+l_2}+\frac{(1-l_1-l_2)^2}{1-y_1/\sqrt{a}-l_1-y_2/\sqrt{a}-l_2}\Bigg)^{1/2}\Bigg)$$
converges to
$${-\frac{1}{2(1-\rho_{12}^2)}\Biggl[\left(\frac{y_1}{\sqrt{\sigma_{11}}}\right)^2+\left(\frac{y_2}{\sqrt{\sigma_{22}}}\right)^2-2\rho_{12} \left(\frac{y_1}{\sqrt{\sigma_{11}}}\right) \left(\frac{y_1}{\sqrt{\sigma_{11}}}\right)\Biggl]}.$$
The last argument follows straightforwardly form the L'Hospital's rule. $\hfill\square$
\endproof

\noindent {\textbf{Proof of Theorem 1:}} Let $(a_n)$ be an
arbitrary sequence such that $a_n \to \infty$.  To simplify the notation, in the argument below, we omit writing the
index $n$ of $a_n$. Assume first that $H(t)=\lambda(t)=t$ (i.e.  $\lambda$ is the Lebesgue measure  on $[0,1]$). Thus  the process (\ref{eq2.188}) reduces to
\begin{eqnarray}
\nonumber D_{\lambda,a}(t)=\sqrt{a}\left(P_{\lambda,a}(t)-t \right).
 \end{eqnarray}
To prove the theorem, we use Lemma \ref{Asy.L2} and Theorem 13.5 of Billingsley (1999). Therefore, we only need to show that for any $0 \le t_1\le t \le t_2 \le 1,$
 \begin{eqnarray}
 \label{tight}  E\left[\left|D_{\lambda,a}(t)-D_{\lambda,a}(t_1)\right|^{2\beta}\left|D_{\lambda,a}(t_2)-D_{\lambda,a}(t)\right|^{2\beta}\right] \le \left|F(t_2)-F(t_1)\right|^{2a},
 \end{eqnarray}
for some  $\beta \ge 0,$ $a>1/2,$ and a nondecreasing continuous function $F$ on $[0,1].$  Take $a=\beta=1$ and $F(t)=t$. We show that
\begin{eqnarray}
  \label{exp1}  E\left[\left(D_{\lambda,a}(t)-D_{\lambda,a}(t_1)\right)^2\left(D_{\lambda,a}(t_2)-D_{\lambda,a}(t)\right)^2\right] \le \frac{ 8a-1}{a^3} \left(t_2-t_1\right)^2.~~~
 \end{eqnarray}
Observe that
  \begin{eqnarray}
  \nonumber  D_{\lambda,a}(t)-D_{\lambda,a}(t_1)=D_{\lambda,a}((t_1,t]) \text { and } D_{\lambda,a}(t_2)-D_{\lambda,a}(t)=D_{\lambda,a}((t,t_2]).
  \end{eqnarray}
  Thus, the expectation in the right hand side of (\ref{exp1}) is equal to
  \begin{eqnarray}
 \label{one} \xi^2(a)E\left[\left\{P_{\lambda,a}((t_1,t])-\lambda((t_1,t])\right\}^2\left\{P_{\lambda,a}((t,t_2])-\lambda((t,t_2])\right\}^2\right],
 \end{eqnarray}
 where $\lambda((t,t_2])=t_2-t$ and $\lambda((t_1,t])=t-t_1$. Expanding the expression
 \begin{eqnarray}
 \label{exp} \left\{P_{\lambda,a}((t_1,t])-\lambda((t_1,t])\right\}^2\left\{P_{\lambda,a}((t,t_2])-\lambda((t,t_2])\right\}^2
  \end{eqnarray}
gives
\begin{align*}
& P^2_{\lambda,a}((t_1,t])P^2_{\lambda,a}((t,t_2])-2\lambda((t,t_2])P^2_{\lambda,a}((t_1,t])P_{\lambda,a}((t,t_2])\\&+
\lambda^2((t,t_2])P^2_{\lambda,a}((t_1,t])-2\lambda((t_1,t])P_{\lambda,a}((t_1,t])P^2_{\lambda,a}((t,t_2])\\&
+4\lambda((t_1,t])\lambda((t,t_2])P_{\lambda,a}((t_1,t])
P_{\lambda,a}((t,t_2])-2\lambda((t_1,t])\lambda^2((t,t_2])P_{\lambda,a}((t_1,t])\\&+\lambda^2((t_1,t])P^2_{\lambda,a}((t,t_2])-2 \lambda^2((t_1,t])\lambda((t,t_2])P_{\lambda,a}((t,t_2])+\lambda^2((t_1,t])\lambda^2((t,t_2]).
\end{align*}
Using the fact that $P_{\lambda,a}(\cdot)$  is a probability measure and $0 \le t_1\le t \le t_2 \le 1$, the expression displayed in  (\ref{exp}) is less than or equal to
\begin{eqnarray*}
&&P^2_{\lambda,a}((t_1,t])P^2_{\lambda,a}((t,t_2])+\lambda^2((t,t_2])
P^2_{\lambda,a}((t_1,t])\\&&+4\lambda((t_1,t])\lambda((t,t_2])P_{\lambda,a}((t_1,t])P_{\lambda,a}((t,t_2])
+\lambda^2((t_1,t])P^2_{\lambda,a}((t,t_2])
\\&&+\lambda^2((t_1,t])\lambda^2((t,t_2])\\&&
\le P_{\lambda,a}((t_1,t])P_{\lambda,a}((t,t_2])+\lambda((t,t_2])
P_{\lambda,a}((t_1,t])
\\&&+4\lambda((t_1,t])\lambda((t,t_2])+\lambda((t_1,t])P_{\lambda,a}((t,t_2])+\lambda((t,t_2])
\lambda((t,t_2]).
\end{eqnarray*}
By (\ref{NIG1}) and (\ref{As.20}) we obtain
\begin{eqnarray}
\nonumber E\left[\left\{P_{\lambda,a}((t_1,t])-\lambda((t_1,t])\right\}^2\left\{P_{\lambda,a}((t,t_2])-\lambda((t,t_2])\right\}^2\right]&\le& \frac{ 8a-1}{a}\times\\&& \label{two} \lambda((t_1,t])\lambda((t,t_2]).~~~~~~~~~~~~~~
\end{eqnarray}
Thus, using (\ref{one}) and (\ref{two}), we have
 \begin{eqnarray*}
 \nonumber  E\left[\left(D_{\lambda,a}(t)-D_{\lambda,a}(t_1)\right)^2\left(D_{\lambda,a}(t_2)-D_{\lambda,a}(t)\right)^2\right]&=&\frac{ 8a-1}{a^3}\lambda(t_1,t]\lambda(t,t_2]\\
\nonumber  &=&\frac{ 8a-1}{a^3}\left(t-t_1\right)\left(t_2-t\right)\\
\nonumber  &\le&\frac{ 8a-1}{a^3} \left(t_2-t_1\right)^2,
 \end{eqnarray*}
for $0 \le t_1\le t \le t_2 \le 1.$ This proves the theorem in the case when $H(t)=t$, i.e. $H$ is the uniform distribution. Observe that, the quantile function $H^{-1}(s)=\inf \left\{t:H(t)\ge s \right\}$ has the property: $H^{-1}(s) \le t$ if and only if $s \le H(t)$. If $U_i$ is uniformly distributed over $[0,1]$, then $H^{-1}(U_i)$ has distribution $H$. Thus, we can use the representation $a_i=H^{-1}(U_i),$ where $(U_i)_{i\ge 1}$  is a sequence of i.i.d. random variables with uniform distribution on $[0,1]$, to have:
\begin{eqnarray}
\nonumber P_{H,a}(t)=P_{\lambda,a}(H(t)) \ \ \text{and} \ \ D_{H,a}(t)=D_{\lambda,a}(H(t))=D_{\lambda,a}\circ H(t), \ \ t\in\mathbb{R},
 \end{eqnarray}
where $P_{\lambda,a}$ is the normalized inverse-Gaussian process with concentration parameter $a$ and Lebesgue base measure $\lambda$ on $[0,1]$. From the uniform case, which was already treated,
we have $D_{\lambda,a}(\cdot)=\sqrt{a}\left(P_{\lambda,a}(\cdot)-\lambda(\cdot)\right)$ $
\overset{d} \to B_\lambda(\cdot) $. Define $\Psi: D[0,1] \to D[-\infty,\infty]$ by $(\Psi x)(t)=x(H(t))$.
Since the function $\Psi$ is uniformly continuous (Billingsley 1999, page 150; Pollard, 1984, page 97), it follows, from the continuous mapping theorem and the fact that $D_{\lambda,a} \overset{d}\to B_\lambda$, that $D_{H,a}=\Psi(D_{\lambda,a}) \overset{d}\to \Psi(B_\lambda)=B_H.$ This completes the proof of the theorem.  $\hfill\square$
 \endproof


\begin{thebibliography}{}

\bibitem{AS} Abramowitz, M., and Stegun, I. (1972). \emph{Handbook of mathematical functions with formulas, graphs, and
mathematical tables}. Dover Publications, Mineola, New York.

\bibitem{BW}  Bickel, P.J. and  Wichura,  M.J. (1971). Convergence criteria for multiparameter stochastic processes and some applications.  \emph{The Annals of Mathematical Statistics}, {42}, 1656-1670.

\bibitem{BF} Bickel, P.J., and Freedman, D.A. (1981). Some asymptotic theory for the bootstrap. \emph{The Annals of Statistics}, {9}, 1196-1217.

\bibitem{B1} Billingsley, P. (1995).  \emph{Probability and Measure}, second edition, John Wiley \& Sons, Inc.

\bibitem{B2} Billingsley, P. (1999). \emph{Convergence of Probability Measures}, third edition. John Wiley \& Sons, Inc.







\bibitem{Das} DasGupta, A. (2008). \emph{Asymptotic Theory of Statistics and Probability}. Springer.




\bibitem{F1}  Ferguson, T.S. (1973). A Bayesian Analysis of Some Nonparametric Problems. \emph{The Annals of Statistics}, {{1}}, 209-230.

\bibitem{FK} Ferguson, T.S., and Klass, M. J. (1972). A Representation of Independent Increment Processes without Gaussian Components. \emph{The Annals of
    Mathematical Statistics}, {1}, 209-230.

\bibitem{IJ} Ishwaran, H., and James, L. F. (2001).  Gibbs Sampling Methods for Stick-Breaking Priors. \emph{Journal of the
American Statistical Association}, {96}, 161-173.

\bibitem{ L.F.} James, L.F (2008). Large sample asymptotics for the two-parameter Poisson-
Dirichlet process. In \emph{Pushing the Lim-
its of Contemporary Statistics: Contributions in Honor of Jayanta K.
Ghosh}, volume 3. Eds. Clarke, B., and Ghosal, S., IMS,  187-199.



\bibitem{K}  Kallenberg, O. (1983). \emph{Random Measures}, third edition. Akademie-Verlag, Berlin.

\bibitem{KB} Kim, N, and Bickel, P. (1987). The limit distribution of a test statistic for bivariate normality. \emph{Statistica Sinica},  {13},
327-349.

\bibitem{LMP} Lijoi, A., Mena, R.H. and Pr\"unster, I. (2005). Hierarchical mixture modelling with normalized
inverse Gaussian priors. \emph{Journal of the American Statistical Association}, {100}, 1278-1291.


\bibitem{F1} Lo , A.Y. (1987). A large sample study of the Bayesian bootstrap. \emph{The Annals of Statistics},  {15},
360-375.




\bibitem{PY} Pitman, J.  and Yor, M.  (1997). The two-parameter Poisson-Dirichlet distribution derived
from a stable subordinator. \emph{The Annals of Probability}, {2}, 855-900.

\bibitem{Po} Pollard, D. (1984). \emph{Convergence of Stochastic Processes}. Springer-Verlag, New York.

\bibitem{RS1}  Resnick, S.I. (1987). \emph{Extreme Values, Regular Variation and Point Processes}. Springer-Verlag, New York.


\bibitem{SE1} Sethuraman, J. and Tiwari, R.C. (1982). Convergence of Dirichlet measures
and the interpretation of their parameter. \emph{Statistical Decision Theory and
Related Topics III}, {2}, 305-315.





%
\bibitem{VW} van der Vaart, A.W., and Wellner, J.A. (1996). \emph{Weak Convergence and Empirical Processes
with Applications to Statistics}. Springer, New York.


\bibitem{ZA} Zarepour, M., and Al Labadi, L. (2012). On a Rapid Simulation of the Dirichlet Process. \emph{Statistics and Probability Letters}, 82, 916-924.


\end{thebibliography}
\end{document}